\documentclass[letterpaper,12pt]{article}
\usepackage{amsmath,amssymb,amsfonts,epsf,epsfig,enumerate,framed}
\usepackage{combelow}

\newtheorem{thm}{Theorem}[section]

\newtheorem{lem}[thm]{Lemma}
\newtheorem{defn}[thm]{Definition}
\newtheorem{alg}[thm]{Algorithm}

\newenvironment{proof}{\noindent {\sc Proof}.}
                {\phantom{a} \hfill \framebox[2.2mm]{ } \bigskip}

\newcommand{\NN}{\mathbb{N}}
\newcommand{\ZZ}{\mathbb{Z}}

\renewcommand{\S}{{\mathcal{S}}}
\renewcommand{\P}{{\mathcal{P}}}
\newcommand{\A}{{\mathcal{A}}}
\newcommand{\T}{{\mathcal{T}}}
\newcommand{\Q}{{\mathcal{Q}}}

\newcommand{\M}{{\mathcal{M}}}
\renewcommand{\phi}{\varphi}
\renewcommand{\P}{{\mathcal{P}}}

\renewcommand{\r}{\hspace*{10mm}}
\newcommand{\dist}{{\rm{dist}}}
\newcommand{\bl}{{\diamond}}

\title{Abstract questionnaires and FS-decision digraphs}
\author{Jiaye Chen, Suzan Kadri, Mateja \v{S}ajna\footnote{Corresponding author. Email: msajna@uottawa.ca. Mailing address: Department of Mathematics and Statistics, University of Ottawa, 150 Louis-Pasteur Private, Ottawa, ON, K1N 6N5, Canada.}, and Ioana \cb{S}chiopu-Kratina \\
{\small University of Ottawa}}

\begin{document}
\maketitle \baselineskip 17pt

\begin{abstract}
A questionnaire is a sequence of multiple choice questions aiming to collect data on a population. We define an {\em abstract questionnaire} as an ordered pair $(N,\M)$, where $N$ is a positive integer  and $\M=(m_0,m_1,\ldots,m_{N-1})$ is an $N$-tuple of positive integers, with  $m_i$, for $i \in \ZZ_N$, as the number of possible answers to question $i$. An abstract questionnaire may be endowed with a skip-list (which tells us which questions to skip based on the sequence of answers to the earlier questions) and a flag-set (which tells us which sequences of answers are of special interest).

An FS-decision tree is a decision tree of an abstract questionnaire that also incorporates the information contained in the skip-list and flag-set. The main objective of this paper is to represent the abstract questionnaire using a directed graph, which we call an FS-decision digraph, that contains the full information of an FS-decision tree, but is in general much more concise. We present an algorithm for constructing a fully reduced FS-decision digraph, and develop the theory that supports it.

In addition, we show how to generate all possible orderings of the questions in an abstract questionnaire that respect a given precedence relation.

\end{abstract}

\section{Introduction}

A questionnaire, from the perspective of this paper, is a sequence of multiple choice questions aiming to collect data on a population. A well-designed questionnaire maximizes the response rate while minimizing the response burden, resulting in the collection of good quality data with minimum bias. The process of designing an effective questionnaire typically relies on its content (that is, the actual questions and possible responses). In this paper, however, the focus is on {\em abstract questionnaires}; that is, we shall not be concerned with the content. Moreover, our main focus will be not the design of the questionnaire itself, but rather on how to represent it by a graph in a way that is as simple and concise as possible, yet retains full information on the questionnaire.

The idea to represent a questionnaire by a graph is not new. In the literature, we find two basic approaches to this task: flow charts (see, for example, \cite{BetHun,Ell,Jab,Par,Pic,KraZamTreMar,StaZin,XinSun}) and decision trees (see \cite{Fen}). A flow chart (called survey chart in \cite{KraZamTreMar}) is a directed graph whose vertices are essentially the questions (each question corresponding to a single vertex, possibly with some additional vertices representing additional controls), and whose arcs (directed edges) represent responses that lead from one question to another. Each arc is labeled with the corresponding response, and is sometimes additionally assigned other parameters (for example, probability of the response), depending on the purpose of the graph model. In a decision tree, however, each question corresponds to a level, and thus may be represented by many vertices, each corresponding to a different sequence of responses to the previous questions. The edges from a given vertex to its children in the tree represent the possible responses to the question. The ordering of the children at each vertex uniquely identifies the responses, so no additional labels on the edges are necessary.

In this paper, we start by modelling a questionnaire with a decision tree, and then condense it into a so-called decision digraph, which retains full information contained in the decision tree. Decision digraphs are a novel way to represent a questionnaire, and can be considered a compromise between decision trees and flow charts: they are potentially much smaller than decision trees, but a lot simpler than flow charts since they do not require any labels on the arcs. Additionally, we endow an abstract questionnaire with a flag-set, which tells us which sequences of answers are of special interest, and a skip-list, which tells us which questions are to be skipped based on the sequence of answers to the earlier questions. As far as we can tell, flagging has not previously been considered in the literature, while skipping has been treated only with the knowledge of the content of the questionnaire (see \cite{FeeFee}). When flagging and skipping are taken into account, we talk of FS-decision trees and FS-decision digraphs. We point out that in both of these representations, the order of the questions is fixed; that is, the order in which the questions are asked does not depend on the sequence of responses to the previous questions. An algorithm for generating all possible orderings of the questions in an abstract questionnaire that satisfies a given precedence relation is discussed in Section~\ref{sec:ordering}.

This paper is organized as follows. In Section~\ref{sec:prerequisites} we give the prerequisites on graphs and relations. After a brief discussion of the ordering of questions in Section~\ref{sec:ordering}, we move on to the main topic, representing the flow of a questionnaire by a graph, presented in Section~\ref{sec:main}. In Subsection~\ref{sec:FS-trees}, we  introduce flag-sets, skip-lists, and FS-decision trees, and present an algorithm for constructing an FS-decision tree of an abstract questionnaire. Vertex equivalence in an FS-decision tree, and FS-decision digraphs are discussed in Subsection~\ref{sec:FS-digraph}; here, we also present an algorithm that constructs an FS-decision digraph, and prove a theorem that shows that the resulting digraph is fully reduced (that is, as small as possible). In the last subsection we talk about how to generate a skip-list and a compatible flag-set that can serve as inputs for our algorithm from a more intuitively constructed pre-skip-list and pre-flag-set. Finally, in Section~\ref{sec:example} we give an example of a simple concrete questionnaire to illustrate the contributions of Sections~\ref{sec:ordering} and \ref{sec:main}.

\section{Prerequisites}\label{sec:prerequisites}

An {\em abstract questionnaire} is an ordered pair $(N,\M)$, where $N$ is a positive integer, and $\M=(m_0,m_1,\ldots,m_{N-1})$ is an $N$-tuple of positive integers. An abstract questionnaire $(N,\M)$ models a questionnaire with $N$ questions (labelled $0,1,\ldots,N-1$) such that question $i$ has $m_i$ possible answers. We denote the set of possible answers to question $i$ by $A_i$; usually, we assume that $A_i=\ZZ_{m_i}$. An abstract questionnaire $(N,\M)$ is called {\em binary} if $\M=(2,2,\ldots,2)$; that is, if every question has exactly two answers.

\subsection{Graphs, trees, and digraphs} \label{sec:graphs}

A {\em graph} $G$ is an ordered triple $(V,E,\psi)$, where $V$ is a non-empty finite set, $E$ is a finite set disjoint from $V$, and $\psi$ is a function assigning to each element of $E$ an unordered pair of elements of $V$. Sets $V$ and $E$ are the {\em vertex set} and {\em edge set}, respectively, of the graph $G$, and $\psi$ is its {\em incidence function}. The elements of $V$ and $E$ are called the {\em vertices} and {\em edges} of $G$, respectively. Vertices $u$ and $v$ in $G$ are said to be {\em adjacent} or {\em neighbours} if $\psi(e)=\{ u,v \}$ for some $e \in E$.

In the definition of a graph, the incidence function may be omitted and edges may be identified with unordered pairs of vertices if there is no ambiguity. We then write $e=\{ u,v \}$, or shortly $e=uv$, instead of $\psi(e)=\{ u,v \}$.

A graph $G'=(V',E',\psi')$ is a {\em subgraph} of a graph $G=(V,E,\psi)$ if $V' \subseteq V$, $E' \subseteq E$, and $\psi'$ is the restriction of $\psi$ to $E'$. If $G'=(V',E',\psi')$ is a subgraph of $G=(V,E,\psi)$  such that $E'=\{  e\in E: \psi(e)=\{u,v\}, u,v \in V' \}$, then $G'$ is the {\em subgraph of $G$ induced by} $V'$, and we write $G'=G[V']$.

A $(v_0,v_k)$-{\em path} (of {\em length} $k$) in a graph $G=(V,E)$ is a sequence $v_0e_1v_1e_2v_2\ldots v_{k-1}e_{k}v_k$, where $v_0,\ldots,v_k$ are distinct vertices of $G$; $e_1,\ldots,e_k$ are edges of $G$, and $e_i=v_{i-1}v_i$ for all $i=1,\ldots,k$.

A {\em cycle} (of {\em length} $k$) in a graph $G=(V,E)$ is a sequence $v_0e_1v_1e_2v_2\ldots v_{k-1}e_{k}v_k$, where $v_0,\ldots,v_{k-1}$ are distinct vertices of $G$ while $v_k=v_0$; $e_1,\ldots,e_k$ are distinct edges of $G$, and $e_i=v_{i-1}v_i$ for all $i=1,\ldots,k$.

A graph $G=(V,E)$ is said to be {\em connected} if for all $u,v \in V$, there exists a $(u,v)$-path in $G$. The {\em distance} $\dist(u,v)$ between vertices $u$ and $v$ in a connected graph $G$ is the length of a shortest $(u,v)$-path in $G$.

A connected graph with no cycles is called a {\em tree}. A {\em rooted tree} is a tree with a distinguished vertex called the {\em root}. If $T$ is a rooted tree with root $r$, and $v$ is any vertex of $T$, then the unique $(r,v)$-path in $T$ defines a sense of direction (away from the root). The vertices of $T$ on this unique $(r,v)$-path, excluding $v$,  are called the ancestors of $v$, and  the neighbours of $v$ not on this path are called the {\em children} of $v$. A vertex is a {\em descendant} of $v$ if $v$ is its ancestor. A rooted tree is said to be {\em ordered} if the set of children at each vertex is ordered.  If $T$ is a rooted tree and $v \in V(T)$, then the {\em subtree of $T$ rooted at $v$} is the subgraph of $T$ induced by the subset of vertices that contains $v$ and all of its descendants.

A {\em directed graph} (shortly {\em digraph}) $D$ is an ordered triple $(V,A,\psi)$, where $V$ is a non-empty finite set, $A$ is a finite set disjoint from $V$, and $\psi$ is a function assigning to each element of $A$ an ordered pair of elements of $V$. Sets $V$ and $A$ are the {\em vertex set} and {\em arc set}, respectively, of the digraph $D$, and $\psi$ is its {\em incidence function}. The elements of $V$ and $A$ are called the {\em vertices} and {\em arcs} of $D$, respectively.
Again, the incidence function may be omitted if there is no ambiguity, and we write $a=(u,v)$  instead of $\psi(a)=(u,v)$.
If $u,v \in V$ are such that $(u,v) \in A$, then $u$ is an {\em in-neighbour} of $v$, and $v$ is an {\em out-neighbour} of $u$. Distinct arcs $a_1,a_2 \in A$ are said to be {\em parallel} if $\psi(a_1)=\psi(a_2)$.

A {\em subdigraph} and {\em vertex-set induced subdigraph} of a digraph are defined analogously to graphs.

A {\em directed cycle} of length $k$ in a digraph $D=(V,A)$ is a sequence $v_0a_1v_1a_2v_2\ldots v_{k-1}a_{k}v_k$, where $v_0,\ldots,v_{k-1}$ are distinct vertices of $D$ while $v_k=v_0$; $a_1,\ldots,a_k$ are arcs of $D$, and $a_i=(v_{i-1},v_i)$ for all $i=1,\ldots,k$. A directed cycle of length one corresponds to an arc of the form $(v,v)$, which is called a {\em directed loop}.

\subsection{Relations, transitive closure, partial orders, and total orders}

A {\em binary relation} on a finite set $S$ is any subset of $S \times S$; that is, a set of ordered pairs of elements from $S$. A binary relation $R$ on a set $S$ is said to be {\em reflexive} if $(a,a) \in R$ for all $a \in S$; {\em irreflexive} if $(a,a) \not\in R$ for all $a \in S$; {\em antisymmetric} if $(a,b),(b,a) \in R$ implies $a=b$, for all $a,b \in S$; and transitive if $(a,b),(b,c) \in R$ implies $(a,c) \in R$, for all $a,b,c \in S$. Elements $a,b \in S$ are said to be {\em comparable} with respect to $R$ if $(a,b) \in R$ or $(b,a) \in R$.

If $R$ and $R'$ are two binary relations on a set $S$, and $R \subseteq R'$, then $R'$ is said to {\em extend} $R$.

The {\em transitive closure} ({\em reflexive closure}, resp.) of a binary relation $R$ on a set $S$ is a minimal binary relation on the set $S$ that is transitive (reflexive, resp.) and extends $R$.

A binary relation that is reflexive, antisymmetric, and transitive is called a {\em partial order}. A {\em minimal element} in a partial order $R$ on the set $S$ is an element $m \in S$ such that $(x,m) \in R$ implies $x=m$, for all $x \in S$.

If $R$ is a partial order on a set $S$, and every pair of elements $a,b \in S$ are comparable with respect to $R$, then $R$ is called a {\em total order}. It is well-known that every partial order has a total order that extends it.

A binary relation $R$ on a set $S$ is equivalent to the digraph $(S,R)$, which has no parallel arcs; conversely, a digraph with no parallel arcs can be viewed as a binary relation. Hence the terms defined in Section~\ref{sec:graphs} for digraphs (for example, directed cycle) apply to binary relations as well.

\section{Ordering questions in a questionnaire} \label{sec:ordering}

In this section, we describe a procedure for constructing all possible orderings of the questions in a questionnaire. Such orderings will be fixed; that is, independent from the responder and independent from the responses to any previous questions. We will, however, assume that these orderings respect a given precedence relation. In other words, our input will be a binary relation $R$ on the set of questions $\ZZ_N$ such that $(a,b) \in R$ if question $a$ must be asked before question $b$. Note that we may assume that $R$ is irreflexive. Our procedure is based on two well-known algorithms: the Roy-Warshall Algorithm for constructing the transitive closure of a binary relation, and Topological Sorting for constructing a total order that extends a given partial order.

The following observation will allow us to determine whether the input precedence relation is admissible (that is, extendible to a total order) or not.

\begin{lem}\label{lem:total}
Let $R$ be an irreflexive binary relation on a set $S$, and $R^*$ its transitive closure. Then the following are equivalent.
\begin{enumerate}[(i)]
\item $R$ has no directed cycles.
\item $R^*$ has no directed loops.
\item $R^*$ has no directed cycles.
\item There exists a total order extending $R$.
\end{enumerate}
\end{lem}

\begin{proof}
(i) $\Rightarrow$ (ii): This is obvious.

(ii) $\Rightarrow$ (iii): This is obvious.

(iii) $\Rightarrow$ (iv): Assume $R^*$ has no directed cycles, and let $\overline{R^*}$ be the reflexive closure of $R^*$. Then $\overline{R^*}$ is reflexive by definition, and transitive since $R^*$ is transitive. We claim $\overline{R^*}$ is antisymmetric. Take any $a,b \in S$ such that $(a,b),(b,a) \in \overline{R^*}$. If $a \ne b$, then $(a,b),(b,a) \in  R^*$, and by transitivity, we have $(a,a) \in R^*$ --- a contradiction since $R^*$ has no directed cycles. Hence the implication $(a,b),(b,a) \in \overline{R^*} \Rightarrow a=b$ holds for all $a,b \in S$.

We conclude that $\overline{R^*}$ is a partial order, and hence there exists a total order $T$ extending $\overline{R^*}$. Since $\overline{R^*}$ extends $R$, we have that $T$ is a total order extending $R$.

(iv) $\Rightarrow$ (i): If $T$ is a total order extending $R$, then $T$ has no directed cycles of length at least 2, and hence $R$ has no directed cycles of length at least 2. Since $R$ is irreflexive, it also has no directed cycles of length 1.
\end{proof}

Note that in all of our algorithms, (ordered) lists are denoted using square brackets, and concatenation of lists is denoted with a plus sign; that is, $[x_1,\ldots,x_n]+[y_1,\ldots,y_m]=[x_1,\ldots,x_n,y_1,\ldots,y_m]$.

\begin{alg}\label{alg:Q}{\rm Constructing all possible orderings of the questions in a questionnaire

\smallskip

\noindent {\bf procedure} Orderings$(N,R)$ \\
\# {\em Input: $N$ is the number of questions, $R$ is an irreflexive binary relation (set of ordered pairs) on $\ZZ_N$.} \\
\# {\em Output: a list $\T$ of all possible total orders on $\ZZ_N$ that extend $R$.} \\
\r $R^*:=$ transitive closure of $R$ \qquad \# {\em use the Roy-Warshall Algorithm} \\
\r {\bf if} $R^*$ has a directed loop {\bf then return} None \\
\r {\bf else} \qquad \# {\em construct the list $\T$ of all total orders} \\
\r\r $R:=$ reflexive closure of $R^*$ \\
\r\r $S:=\ZZ_N$ \qquad \# {\em the underlying set} \\
\r\r $L:=[ \; ]$ \qquad \# {\em current total order} \\
\r\r $\T:=[ \; ]$ \\
\r\r TotalOrder$(S,R,L,\T)$ \\
\r\r {\bf return} $\T$

\smallskip

\noindent {\bf procedure} TotalOrder$(S,R,L,\T)$ \\
\# {\em Input: $S$ is the underlying set, $R$ is a partial order (set of ordered pairs) on $S$, $L$ is the current total order we are constructing, $\T$ is the list of total orders constructed so far.} \\
\# {\em Output: the list $\T$ of all possible total orders on $S$ that extend $R$.} \\
\r {\bf if } $S=\emptyset$ {\bf then} $\T=\T + [L]$ \qquad \# {\em total order $L$ is complete; add it to list $\T$} \\
\r {\bf else} \\
\r\r $M:=$ the set of minimal elements of the partial order $(S,R)$ \\
\r\r {\bf for} $m \in M$ {\bf do} \\
\r\r\r $L':=L + [m]$ \\
\r\r\r $S':=S-\{ m \}$ \\
\r\r\r $R':=R \cap (S' \times S')$ \\
\r\r\r TotalOrder$(S',R',L',\T)$  \\
}
\end{alg}


\section{Representing the flow of a questionnaire by a graph}\label{sec:main}

In this section, we aim to represent the flow of a questionnaire by a graph. We shall start with the well-known decision tree representation, then introduce the FS-decision tree, that is,  a more compact version of the decision tree that also encodes information on special answer sequences, and finally describe an FS-decision digraph, a much more compact graph that nevertheless contains the full information on the questionnaire. Throughout this section, we shall assume that the order of the questions has been fixed; that is, question 0 is asked first, then question 1, and so on until question $N-1$. An additional question $N$, with no answers, will represent the end.

Throughout this section, we shall assume that we have an  abstract questionnaire $\Q=(N,\M)$, and that $A_i=\ZZ_{m_i}$ the set of possible answers to question $i$. In addition, for all $i \in \ZZ_N$, we define $A_i^*=A_i \cup \{ * \}$. The symbol $*$ will represent all possible answers, and hence may correspond to a skipped question.

\begin{defn}{\rm
Let $k$ and $\ell$ be integers, $0 \le k \le \ell \le N-1$.
A {\em $(k,\ell)$-answer string} for an  abstract questionnaire $\Q=(N,\M)$ is an element of $A_k^* \times A_{k+1}^* \times \ldots \times A_{\ell}^*$; that is, a string of the form $a_k a_{k+1} \ldots a_{\ell}$, where $a_i \in A_i^*$ for all $i=k,k+1,\ldots,\ell$.

An {\em answer string} for $\Q$ is either the empty string, denoted $\epsilon$, or a $(k,\ell)$-answer string for $0 \le k \le \ell \le N-1$.
}
\end{defn}


The concatenation of strings $a$ and $b$ will be denoted $ab$. If $a$ is a $(k,\ell)$-answer string and $b$ is an $(\ell+1,\ell')$-answer string, for some $0 \le k \le \ell < \ell' \le N-1$, then $ab$ is a $(k,\ell')$-answer string.

We now formally define a decision tree for an abstract questionnaire $\Q=(N,\M)$, as well as an assignment of questions and answer strings to its vertices.

\begin{defn}{\rm
A {\em decision tree} $T$ for the abstract questionnaire $\Q=(N,\M)$ is an ordered rooted tree with $V(T)=\ZZ_n$, where $n=1+\sum_{i=0}^{N-1} m_0 m_1 \ldots m_{i}$.
We define the edge set of $T$, together with the {\em question assignment} $\kappa: V(T) \to \ZZ_{N+1}$ and {\em answer string assignment} $\alpha: V(T) \to \{ \epsilon \} \cup \bigcup_{i=0}^{N-1} A_0 \times A_1 \times \ldots \times A_i$ recursively as follows.
\begin{enumerate}[(i)]
\item Vertex 0 is the root of $T$, and $\kappa(0)=0$ and $\alpha(0)=\epsilon$.
\item A vertex $u$ with $\kappa(u)=q$ has $m_{q}$ children if $q < N$, and has no children if $q=N$.
\item If $u$ is a vertex with $\kappa(u)=q$ and $c$ is the $j$-th child of vertex $u$ (for $j \in \ZZ_{m_q}$),  then $\kappa(c)=\kappa(u)+1$ and $\alpha(c)=\alpha(u)j$.
\end{enumerate}
}
\end{defn}

Observe that a vertex $u$ with $\kappa(u)=q$ is at distance $q$ from the root, and its answer string $\alpha(u)$ is a $(0,q-1)$-answer string.

The decision tree of an abstract questionnaire is easily constructed using a Breadth-First-Search algorithm; see Algorithm~\ref{alg:S-dt} with $F=\emptyset$ and $\S=(\emptyset, \ldots, \emptyset)$.

\subsection{FS-decision trees} \label{sec:FS-trees}

In this section, we shall upgrade our decision tree for an abstract questionnaire in two ways. First, we are going to skip all unnecessary vertices; these are vertices corresponding to answer strings that we do not want to be included in the questionnaire.  Second, we are going to flag all those vertices whose answer strings are of special interest (possibly contradictory, but should not be excluded). The latter feature, of course, has no effect on the structure of the tree. We call the resulting tree an FS-decision tree of the abstract questionnaire.

To that end, we make two additional assumptions on the input. First, we assume that we have, for each question $q \ge 1$, a set $S_q$ of $(0,q-1)$-answer strings.  The answer strings in $S_q$ represent sequences of answers to questions 0 to $q-1$ for which question $q$ should be skipped. The list $(S_0,S_1,\ldots,S_{N})$ will be called a skip-list (see Definition~\ref{def:skip-list} below).
Second,  we assume that we are given a set of answer strings that are of special interest (contradictory or in any other way significant so that we wish to keep track of them); we call this set a flag-set (see Definition~\ref{def:flag-set} below). We will then flag every vertex of the tree whose answer string is of special interest.  The restrictions imposed on a skip-list in Definition~\ref{def:skip-list} will ensure that an FS-decision tree has a unique skip-list, and the additional assumption of compatibility that we impose on a flag-set later on (see Definition~\ref{def:compatibleF}) will guarantee uniqueness of the flag-set. Both of these properties will consequently improve the efficiency of our algorithm for constructing an FS-decision digraph in Section~\ref{sec:FS-digraph}.

We formally define a flag-set and a skip-list as follows.

\begin{defn}\label{def:flag-set}{\rm
A {\em flag-set} for the abstract questionnaire $\Q=(N,\M)$ is a set $F$ of $(0,N-1)$-answer strings.
A $(0,k)$-answer string $a_0 \ldots a_k$ of $\Q$ is said to be {\em flagged} with respect to the flag-set $F$ if  for some $i$, $0 \le i \le k$, the $(0,N-1)$-answer string $a_0 \ldots a_i*\ldots*$ is in $F$.
}
\end{defn}

\begin{defn}\label{def:skip-list}{\rm
A {\em skip-list} for the abstract questionnaire $\Q=(N,\M)$ is an $(N+1)$-tuple $\S=(S_0,S_1,\ldots,S_{N})$ such that the following hold.
\begin{enumerate}[(i)]
\item $S_0 =\emptyset=S_N$.
\item For each $q$, $1 \le q \le N$, we have that $S_q$ is a set of $(0,q-1)$-answer strings.
\item For each $q$, $1 \le q \le N$, and each $(0,q-1)$-answer string $a$:
\begin{itemize}
\item if $a \in S_{q}$, then for all $j \in A_{m_q}$,  we have that $aj \not\in S_{q+1}$ and for all $b \in A_{m_{q+1}}^* \times \ldots \times A_{m_{k}}^*$, $q+1 \le k \le N-1$,  the $(0,k)$-answer string $aj b$ is not in $S_{k+1}$; and
 \item if $a \not\in S_{q}$, then  $a* \not\in S_{q+1}$ and for all $b \in A_{m_{q+1}}^* \times \ldots \times A_{m_{k}}^*$, $q+1 \le k \le N-1$,
 the $(0,k)$-answer string $a* b$ is not in $S_{k+1}$.
 \end{itemize}
\end{enumerate}
}
\end{defn}

Note that Requirement (iii) in Definition~\ref{def:skip-list} will ensure that the set $S_q$ contains precisely those $(0,q-1)$-answer strings that correspond to the skipped vertices of the FS-decision tree (see Definition~\ref{def:S-dt} below).

We are now ready to define an FS-decision tree.

\begin{defn}\label{def:S-dt}{\rm
An {\em FS-decision tree} $T$ for the abstract questionnaire $\Q=(N,\M)$ with a flag-set $F$ and skip-list $\S=(S_0,\ldots,S_{N})$ is an ordered rooted tree with $V(T)=\ZZ_n$, for some integer $n \le 1+\sum_{i=0}^{N-1} m_0 m_1 \ldots m_{i}$.
We define the edge set of $T$, together with the
\begin{itemize}
\item subset $U$ of {\em  skipped vertices},
\item {\em question assignment} $\kappa: V(T) \to \ZZ_{N+1}$ and
\item {\em answer string assignment} $\alpha: V(T) \to \{ \epsilon \} \cup \bigcup_{i=0}^{N-1} A_0^* \times A_1^* \times \ldots \times A_i^*$
\end{itemize}
recursively as follows.
\begin{enumerate}[(i)]
\item Vertex 0 is the root, $\kappa(0)=0$ and $\alpha(0)=\epsilon$.
\item If $u$ is a vertex with $\kappa(u)=q \le N$, then
    $u \in U$ if and only if $\alpha(u) \in S_{q}$.
\item If $u$ is a vertex with $\kappa(u)=q<N$, then
\begin{itemize}
\item if $u \in U$, then $u$ has exactly one child, say $c$, and $\kappa(c)=q+1$ and $\alpha(c)=\alpha(u)*$;
\item if $u \not\in U$, then $u$ has exactly $m_q$ children, and for each $j \in A_q$, the $j$-th child, say $c$, satisfies $\kappa(c)=q+1$ and $\alpha(c)=\alpha(u)j$.
\end{itemize}
\item A vertex $u$ with $\kappa(u)=N$ has no children.
\end{enumerate}
In addition, we define the {\em flag function} $\phi: V(T) \to \ZZ_2$ as follows: for all $v \in V(T)$, we have $\phi(v)=1$ (and we say that $v$ is {\em flagged} with respect to $F$) if and only if $\alpha(v)$ is flagged with respect to $F$.
}
\end{defn}

From the above definition, it is clear that the skip-list of an abstract questionnaire $\Q=(N,\M)$ uniquely defines its FS-decision tree.

In Lemma~\ref{lem:uniqueS} below, we show the converse to  this statement: namely, that each FS-decision tree arises from a unique skip-list. In the proof, we shall use the following lemma.

\begin{lem}\label{lem:skip-list-and-vx}
Let $\Q=(N,\M)$ be an abstract questionnaire with a skip-list $\S=(S_0,\ldots,S_{N})$, and $T$ its FS-decision tree with answer string assignment $\alpha$. For $0 \le k \le N-1$, let $a$ be a $(0,k)$-answer string such that $a \in S_{k+1}$. Then there exists $u \in V(T)$ such that $\alpha(u)=a$.
\end{lem}

\begin{proof}
Let $a=a_0 \ldots a_k$, and let $i$, for $0 \le i \le k$, be the largest index such that there exists $v \in V(T)$ with $\alpha(v)=a_0 \ldots a_i$. Suppose $i<k$.
Then there is no vertex in $T$ with answer string $a_0 \ldots a_{i} a_{i+1}$. Suppose first that $a_{i+1} \ne *$. Then, by Definition~\ref{def:S-dt}, we know that $v \in U$ and $\alpha(v) \in S_{i+1}$, and  by requirement (iii) of Definition~\ref{def:skip-list}, we have that $a=\alpha(v) a_{i+1} \ldots a_k \not\in S_{k+1}$, which is a contradiction. Hence $a_{i+1}=*$, and since there is no vertex in $T$ with answer string $a_0 \ldots a_{i} a_{i+1}$,  Definition~\ref{def:S-dt} tells us that $v \not\in U$. Hence $\alpha(v) \not\in S_{i+1}$, and it follows that $a= \alpha(v) * a_{i+2} \ldots a_k \not\in S_{k+1}$, again a contradiction.

We conclude that $i=k$, and hence there exists $u \in V(T)$ such that $\alpha(u)=a$.
\end{proof}

\begin{lem}\label{lem:uniqueS}
Let $\Q=(N,\M)$ be an abstract questionnaire with a skip-list $\S=(S_0,\ldots,S_{N})$, and $T$ its FS-decision tree. If $T$ is also an FS-decision tree for $\Q$ with skip-list $\S'=(S_0',\ldots,S_{N}')$, then $\S=\S'$.
\end{lem}

\begin{proof}
Let $U$, $\kappa$, and $\alpha$ be the set of skipped vertices, question assignment, and answer string assignment, respectively,  for the FS-decision tree $T$. Observe that, by Definition~\ref{def:S-dt}, a vertex $v \in V(T)$ is in $U$ if and only if $\alpha(v) \in S_{\kappa(v)+1}$, and also   if and only if $\alpha(v) \in S'_{\kappa(v)+1}$.

Clearly, $S_0=S_0'$ and $S_N=S_N'$ by definition. We show that $S_{k+1}=S'_{k+1}$ for all $k=0, \ldots, N-2$.
Take any $k \in \{ 0, \ldots, N-2 \}$, and suppose $S_{k+1}' \ne S_{k+1}$. Then, without loss of generality, there exists a $(0,k)$-answer string $a=a_0 \ldots a_k$ such that $a \in S_{k+1}'-S_{k+1}$.

 Since $a \in S'_{k+1}$, by Lemma~\ref{lem:skip-list-and-vx}, there is a vertex $u \in V(T)$ such that $\alpha(u) = a$. However, by the first observation of this proof, we know $u \in U$, and hence $a \in S_{k+1}$, a contradiction.

We conclude that $S_{k+1}' = S_{k+1}$.
\end{proof}

From Definition~\ref{def:S-dt}, it is also clear that the flag function of an S-decision tree is uniquely determined by the flag-set. For the converse to hold, an additional condition on the flag-set --- namely, compatibility with the skip-list, to be defined below --- is required.

\begin{defn}\label{def:compatibleF}{\rm
Let $F$ be a flag-set of an abstract questionnaire $\Q=(N,\M)$ with skip-list $\S=(S_0,\ldots,S_{N})$. Then $F$ is said to be {\em compatible with $\S$} if for each $f \in F$ there exists $\ell \in \{ 0, \ldots, N-1 \}$ such that the following hold:
\begin{enumerate}[(i)]
\item $f=f_0 \ldots f_{\ell} * \ldots *$;
\item for all $i \in \{ 0, \ldots, \ell \}$, we have
$$f_0 \ldots f_{i-1} \in S_i \Longleftrightarrow f_i=*;$$
\item if $f_{\ell} \ne *$, then for some $j \in A_{\ell}$, the $(0,N-1)$-answer string $f_0 \ldots f_{\ell-1}j * \ldots *$ is not in $F$; and
\item for all $i \in \{ 0, \ldots, \ell-1\}$, we have that the $(0,N-1)$-answer string $f_0 \ldots f_i * \ldots *$ is not in $F$.
\end{enumerate}
}
\end{defn}
Note that for $i=0$, Statement (ii) says that $\epsilon \in S_0 \Longleftrightarrow f_0=*$, in other words, that $f_0 \ne *$.

The following lemma will explain the significance of the parameter $\ell$ in Definition~\ref{def:compatibleF}.

\begin{lem}\label{lem:flagged}
Let $F$ be a flag-set and $\S=(S_0,\ldots,S_{N})$ a skip-list of an abstract questionnaire $\Q=(N,\M)$, and assume that $F$ is compatible with $\S$. Take any $(0, N-1)$-answer string $f_0 \ldots f_{\ell} * \ldots * \in F$, where $\ell$ satisfies the conditions in Definition~\ref{def:compatibleF}. Then the following hold.
\begin{enumerate}[(a)]
\item For all $i \in \{ 0, \ldots, \ell \}$, there exists $v \in V(T)$ such that $\alpha(v)=f_0 \ldots f_i$.
\item Let $v \in V(T)$ be such that $\alpha(v)=f_0 \ldots f_{\ell}$. Then $v$ is flagged, every descendant of $v$ is flagged, and every ancestor of $v$ is not flagged.
\end{enumerate}
\end{lem}

\begin{proof}
\begin{enumerate}[(a)]
\item We prove this statement by induction on $i$. Since $\epsilon \not\in S_0$, by Definition~\ref{def:S-dt}, the root vertex has $m_0$ children, with answer strings $0, \ldots, m_0-1$. Furthermore, by Definition~\ref{def:compatibleF}, we know that $f_0 \ne *$. Hence indeed, there is a vertex with  answer string $f_0$.

Suppose the claim holds for some $i$, $0 \le i <\ell$. Let $v \in V(T)$ be such that $\alpha(v)=f_0 \ldots f_i$. By Definition~\ref{def:S-dt}, if $f_0 \ldots f_i \in S_{i+1}$, then $v$ has a child with answer string $f_0 \ldots f_i *$, and by Definition~\ref{def:compatibleF}, we know that $f_{i+1} = *$. If, however, $f_0 \ldots f_i \not\in S_{i+1}$, then by Definition~\ref{def:S-dt}, vertex $v$ has a child with answer string $f_0 \ldots f_i j$ for all $j \in A_{i+1}$, and by Definition~\ref{def:compatibleF}, we know that $f_{i+1} \ne *$. In both cases, we conclude that $v$ has a child with answer string $f_0 \ldots f_i f_{i+1}$.

The claim then follows by induction.

\item Note that $v$ exists by Statement (a). Let $u$ be any ancestor of $v$. Then $\alpha(u)=f_0 \ldots f_{i}$ for some $i<\ell$, and by Definition~\ref{def:compatibleF}, we know that the $(0,N-1)$-answer string $f_0 \ldots f_{j}* \ldots *$ is not in $F$ for all $j < \ell$. Hence by Definition~\ref{def:flag-set}, vertex $u$ is not flagged.

Let $u$ be a vertex of $T$ that is either $v$ itself or a descendant of $v$. Then $\alpha(u)=f_0 \ldots f_{\ell} f_{\ell+1} \ldots f_{k}$ for some $k$, $\ell \le k \le N-1$, and $f_{\ell+1}\ldots f_{k} \in A_{\ell+1}^* \times \ldots \times A_{k}^*$. Then by Definition~\ref{def:flag-set}, since $f_0 \ldots f_{\ell} * \ldots * \in F$, we know that $u$ is flagged.
\end{enumerate}
\end{proof}

We are now ready to outline our algorithm for constructing an FS-decision tree of a an abstract questionnaire.


\begin{alg}\label{alg:S-dt}{\rm Constructing an FS-decision tree of an abstract questionnaire

\smallskip

\noindent {\bf procedure} FS-tree$(N,\M,F,\S)$ \\
\# {\em Input: abstract questionnaire $(N,\M)$ with flag-set $F$ and skip-list $\S=(S_0, \ldots, S_N)$.}\\
\# {\em Output: FS-decision tree $T$ for $(N,\M)$ with the subset of skipped vertices $U$, question assignment $\kappa$, answer string assignment $\alpha$, and flag function $\phi$.} \\
\# {\em The vertices of $T$ are labelled $0,1,2,\ldots$.} \\
\# {\em {\rm Out}$(u)$ is the list of children of vertex $u$.} \\
\r $U:= \emptyset$ \\
\r $\kappa(0):=0$ \\
\r $\alpha(0):=\epsilon$ \\
\r $\phi(0):=0$  \\
\r $L:=[0]$  \qquad \# {\em BFS queue of unprocessed vertices} \\
\r $c:=0$  \qquad \# {\em last vertex label used}\\
\r {\bf while} $L \ne [\;]$ {\bf do} \\
\r\r $u:=$ first vertex in $L$ \\
\r\r remove $u$ from $L$ \\
\r\r $q:=\kappa(u)$ \\
\r\r Out$(u):=[\;]$ \\
\r\r {\bf if} $q <N$ {\bf then} \\
\r\r\r {\bf if} $u \in U$ {\bf then} $A=\{ * \}$\\
\r\r\r {\bf else} $A=\ZZ_{m_q}$\\
\r\r\r {\bf for all} $a_q \in A$ {\bf do} \\
\r\r\r\r $c:=c+1$ \\
\r\r\r\r $\kappa(c)=q+1$ \\
\r\r\r\r $\alpha(c)=\alpha(u)a_q$ \\
\r\r\r\r Out$(u):=$Out$(u) + [ c ]$ \\
\r\r\r\r {\bf if} $\alpha(c)\in S_{\kappa(c)}$  {\bf then} $U := U \cup \{ c \}$ \qquad \# {\em $c$ is a skipped vertex} \\
\r\r\r\r {\bf if } flagged$(\alpha(c),F)=1$ {\bf then} $\phi(c):=1$ \\
\r\r\r\r {\bf else } $\phi(c):=0$  \\
\r\r\r\r $L:=L + [c]$ \\
\r {\bf return} $c$, Out, $U$, $\kappa$, $\alpha$, $\phi$ \qquad \#  $V(T)=\{ 0,1,\ldots,c \}$

\bigskip

\noindent {\bf procedure} flagged$(a_0 \ldots a_k,F)$ \\
\# {\em Input: a $(0,k)$-answer string $a_0 \ldots a_k$, flag-set $F$.}\\
\# {\em Output: 1 if $a_0 \ldots a_k$ is flagged; 0 otherwise.} \\
\r Ans$:=0$ \\
\r {\bf for } $i:=0$ {\bf to }$k$  \\
\r\r {\bf if } the $(0,N-1)$-answer string $a_0 \ldots a_i * \ldots * \in F$ \\
\r\r {\bf then } \\
\r\r\r Ans$:=1$  \\
\r\r\r {\bf break} \\
\r {\bf return} Ans
}
\end{alg}

\subsection{FS-decision digraphs}\label{sec:FS-digraph}

It may happen --- and very frequently it does --- that one subtree of an FS-decision tree looks just like another subtree, except for the initial substring of all corresponding answer strings. In this case, it is advantageous to merge these two subtrees and thus reduce the size of the graph. There are two things to be careful about. First, we do not wish to lose any information; to that end, we replace the answer string of each vertex with the set of answer strings of all vertices it represents. Second, the resulting graph will no longer be a tree. However, if we imagine our FS-decision tree as a digraph with all edges directed away from the root, then our merged graph retains this sense of direction away from the root; even though it may not be acyclic, it does not have any directed cycles. (In fact, every cycle decomposes into two paths of the same length directed in opposite ways.) This sense of direction is, of course, the flow of the questionnaire.

\subsubsection{Vertex equivalence and FS-decision digraphs}

The following definition will make the idea of similar subtrees precise.

\begin{defn}\label{def:eq-vx}{\em
Let $\Q=(N,\M)$ be an abstract questionnaire with flag-set $F$ and skip-list $\S$, and let $T$ be the FS-decision tree for $\Q$, with question assignment $\kappa$ and answer string assignment $\alpha$. Let $v$ and $w$ be two vertices of $T$, and $T_v$ and $T_w$ the subtrees of $T$ rooted at $v$ and $w$, respectively. Furthermore, assume $\kappa(v)=\kappa(w)=k$.

Then vertices $v$ and $w$ are said to be {\em equivalent} in $T$ if there exists a bijection $\Phi: V(T_v) \to V(T_w)$ such that for each $u \in V(T_v)$,
\begin{enumerate}[(a)]
\item if $\alpha(u)=\alpha(v)x$, where $x=\epsilon$ or $x$ is a $(k,\kappa(u)-1)$-answer string, then $\alpha(\Phi(u))=\alpha(w)x$; and
\item $\alpha(u)$ is flagged $\Longleftrightarrow \alpha(\Phi(u))$ is flagged.
\end{enumerate}
}
\end{defn}

The following properties of the mapping $\Phi$ from Definition~\ref{def:eq-vx} are easy to see.

\begin{lem}\label{lem:phi}
Let $\Q=(N,\M)$ be an abstract questionnaire with flag-set $F$ and skip-list $\S$, and let $T$ be the FS-decision tree for $\Q$, with the set of skipped vertices $U$, question assignment $\kappa$, and answer string assignment $\alpha$. Let $v$ and $w$ be two vertices of $T$, and $T_v$ and $T_w$ the subtrees of $T$ rooted at $v$ and $w$, respectively. Furthermore, let $\kappa(v)=\kappa(w)=k$, and let  $\Phi: V(T_v) \to V(T_w)$ be a bijection satisfying Property (a) from Definition~\ref{def:eq-vx}.  Then:
\begin{enumerate}[(i)]
\item $\Phi(v)=w$;
\item for all $u \in V(T_v)$, we have $\kappa(\Phi(u))=\kappa(u)$;
\item $\Phi$ is unique; 
\item $\Phi$ is an isomorphism from $T_v$ to $T_w$; and
\item for all $u \in V(T_v)$, if $m_{\kappa(u)} \ge 2$, then   $u \in U$ if and only if $\Phi(u) \in U$.
\end{enumerate}
\end{lem}

\begin{proof}
\begin{enumerate}[(i)]
\item Since $T$ has a unique vertex with answer string $\alpha(w)$, this statement follows from Property (a) in Definition~\ref{def:eq-vx}, with $x=\epsilon$.
\item Let $u \in V(T_v)$. Since, by assumption, $\kappa(v)=\kappa(w)$, by Property (a),  the answer strings of $u$ and $\Phi(u)$ are of the same length. Hence $\kappa(u)=\kappa(\Phi(u))$.
\item Since no two distinct vertices in an FS-decision tree have the same answer string, this statement follows directly from Property (a).
\item By (i) and Property (a), the bijection $\Phi$ maps $v$ to $w$, and for any $u \in V(T_v)$, it maps the children of vertex $u$ to the children of vertex $\Phi(u)$. Hence $\Phi$ is an isomorphism from $T_v$ to $T_w$.
\item  Take any $u \in V(T_v)$, let $\ell=\kappa(u)$, and assume $m_{\ell} \ge 2$.
    If  $u \in U$, then $\alpha(u) \in S_{\ell}$ and $\ell \le N-1$. By Definition~\ref{def:S-dt}, vertex $u$ has exactly one child, and since $\Phi$ is an isomorphism from $T_v$ to $T_w$ that maps the root of $T_v$ to the root of $T_w$, it also maps the children of $u$ to the children of $\Phi(u)$. It follows that $\Phi(u)$ has exactly one child, and since $\kappa(\Phi(u))=\ell$ and $m_{\ell} \ge 2$, it follows that  $\Phi(u) \in U$.

    Since $\Phi^{-1}$ is an isomorphism from $T_w$ to $T_v$ that satisfies Property (a) from Definition~\ref{def:eq-vx}, the converse follows by symmetry.
\end{enumerate}
\end{proof}

In other words, vertices $v$ and $w$ of an FS-decision tree are equivalent if there exists an isomorphism from $T_v$ to $T_w$ that preserves the flagging property as well as preserves each answer string minus the initial substring $\alpha(v)$.

A digraph resulting from merging some of the pairs of equivalent vertices will be called an FS-decision digraph of the abstract questionnaire; see Definition~\ref{def:S-D} below. Any FS-decision digraph is an ordered levelled digraph, in the following sense.

\begin{defn}\label{def:levD}{\em
An {\em ordered levelled digraph} is a digraph $D=(V,A)$ together with
\begin{enumerate}[(a)]
\item an ordering of the set of out-neighbours of each vertex, and
\item a {\em level function} $\kappa: V \to \NN$ assigning to each vertex $v$ its level $\kappa(v)$ so that
\begin{enumerate}[(i)]
\item there exists at least one vertex $v$ such that $\kappa(v)=0$, and
\item if $(v,u) \in A$, then $\kappa(u)=\kappa(v)+1$.
\end{enumerate}
\end{enumerate}
}
\end{defn}

Note that an ordered rooted tree can be viewed as an ordered levelled digraph with $\kappa$ being the distance from the root, all arcs being directed away from the root, and each set of out-neighbours ordered as the corresponding set of children in the tree. The analogue of a subtree rooted at a vertex $v$ in a rooted tree is the subdigraph rooted at $v$, defined below. It is easy to see that, as stated in the subsequent lemma, a subdigraph rooted at $v$ of an ordered levelled digraph is itself an ordered levelled digraph.

\begin{defn}\label{def:subD}{\em
Let $D=(V,A)$ be an ordered levelled digraph with  level function $\kappa$, and $v \in V$. The {\em subdigraph of $D$ rooted at $v$} is the induced digraph $D_v=D[V']$ where $V'=\{ u \in V: \mbox{there exists a directed } (v,u)-\mbox{path in } D\}$.
}
\end{defn}

\begin{lem}
Let $D=(V,A)$ be an ordered levelled digraph with  level function $\kappa$, and $v \in V$. Using the same ordering of the set of out-neighbours at each vertex, the subdigraph of $D$ rooted at $v$ is an ordered levelled digraph  with level function $\kappa_v=\kappa-\kappa(v)$.
\end{lem}

We are now ready to define an FS-decision digraph of an abstract questionnaire. Note that the arc set of an ordered levelled digraph, as well as the ordering of the out-neighbours of each vertex, is fully defined by giving an (ordered) list of out-neighbours, to be denoted Out$(v)$, for each vertex $v$. Note that  Out$(v)$ may contain repeated elements.

\begin{defn}\label{def:S-D}{\em
Let $\Q=(N,\M)$ be an abstract questionnaire with flag-set $F$ and skip-list $\S$, and let $T$ be the FS-decision tree for $\Q$, with set of skipped vertices $U$, question assignment $\kappa$, answer string assignment $\alpha$, and flag function $\phi$.

The {\em FS-decision digraph $D_T$ corresponding to $T$} is the orientation of $T$ with all arcs directed away from the root. More precisely, $D_T$ is an ordered levelled digraph with level function $\kappa$, and the ordering of the set of out-neighbours of each vertex inherited from $T$. The {\em answer-string assignment} $\A$ for $D_T$ is  defined as $\A(v)=\{ \alpha(v) \}$ for all $v \in V(T)$. The set of skipped vertices and flag function are inherited from $T$.

An {\em FS-decision digraph} for $\Q$ is defined recursively as follows: it is either the FS-decision digraph $D_T$ corresponding to $T$, or is obtained from any FS-decision digraph $D$ via the operation merge$(v,w)$, where $v$ and $w$ are vertices of $D$ equivalent in $T$. The operation merge$(v,w)$ is defined as follows:
\begin{enumerate}[(a)]
\item delete the subdigraph $D_w$ of $D$ rooted at $w$;
\item for each in-neighbour $p$ of $w$, replace each occurrence of $w$ in Out$(p)$ with $v$;
\item for each $u \in V(D_v)$, if $z$ is the corresponding vertex in $D_w$ (that is, $z=\Phi(u)$ where $\Phi: T_v \to T_w$ is the unique isomorphism satisfying Definition~\ref{def:eq-vx}), then adjoin the elements of $\A(z)$ to $\A(u)$;
\item delete all $z \in V(D_w)$ from $U$;
\item restrict $\kappa$ and $\phi$ to $V(D)-V(D_w)$.
\end{enumerate}

An FS-decision digraph for $\Q$ is said to be {\em fully reduced} if it contains no pair of distinct vertices that are equivalent in $T$.
}
\end{defn}

In other words, $D_T$ is the FS-decision digraph obtained from $T$ in the obvious way, that is, by orienting each edge away from the root. Any other FS-decision digraph is obtained from $D_T$ via a sequence of operations whereby the subdigraphs rooted at two equivalent vertices are merged while preserving all information encoded by the FS-decision tree. Recall that the list of out-neighbours of a vertex in an FS-decision digraph may contain repeated elements; in other words, an FS-decision digraph may contain parallel arcs.

\begin{figure}[t!]
\centerline{\includegraphics[scale=0.7]{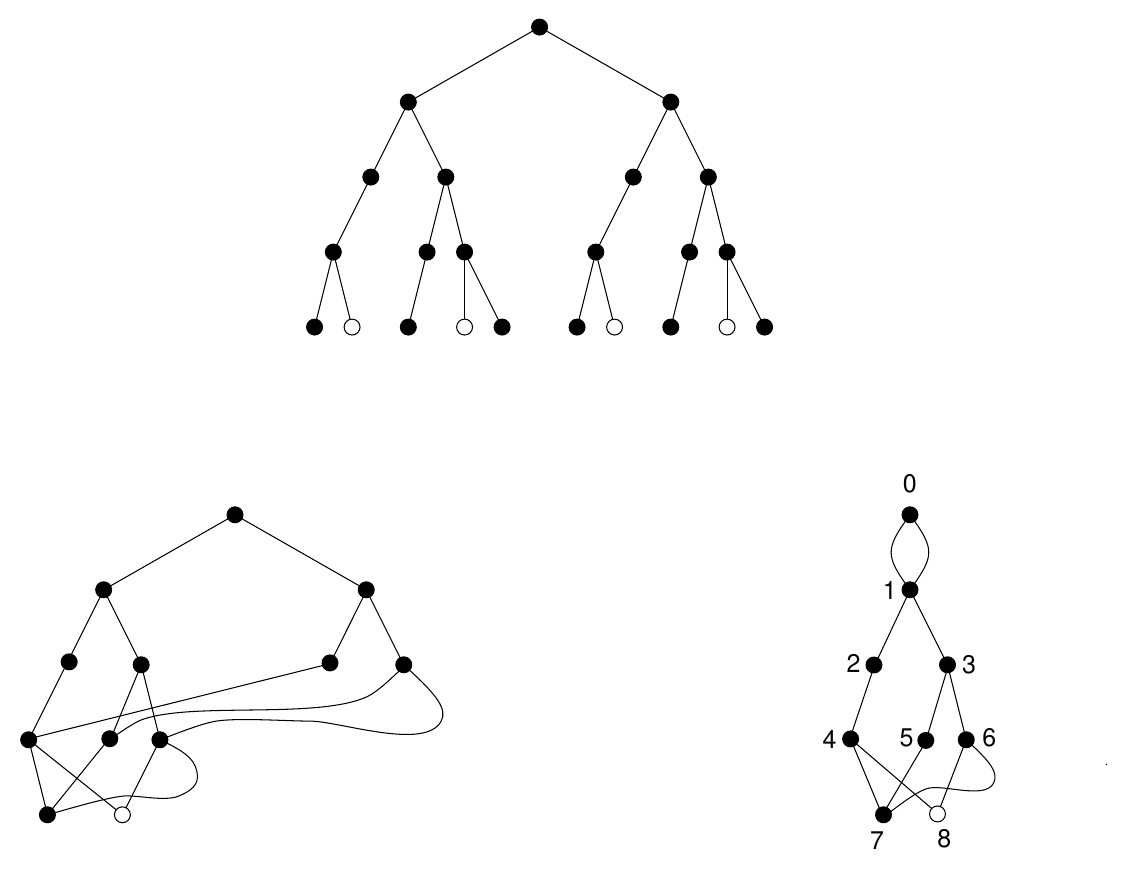}}
\caption{The FS-decision tree $T$ (top) for a binary questionnaire $\Q$. If we imagine all edges directed downwards, this figure also represents the FS-decision digraph $D_T$ corresponding to $T$.
Bottom left: an FS-decision digraph $D'$ for $\Q$ that is not fully reduced. Bottom right: the fully reduced FS-decision digraph $D''$ for $\Q$ . All edges are directed downwards, out-neighbours are ordered from left to right, and flagged and unflagged vertices are coloured white and black, respectively. Vertices in $D''$ are labelled in the order created by Algorithm~\ref{alg:S-dd}.}\label{fig:eg}
\end{figure}

Figure~\ref{fig:eg} shows, for an example of a binary questionnaire $\Q$, its FS-decision tree $T$ and the FS-decision digraph corresponding to $T$ (both represented by the tree on top), an FS-decision digraph for $\Q$ that is not fully reduced, and the fully reduced FS-decision digraph for $\Q$.

\begin{figure}[t!]
\centerline{\includegraphics[scale=0.7]{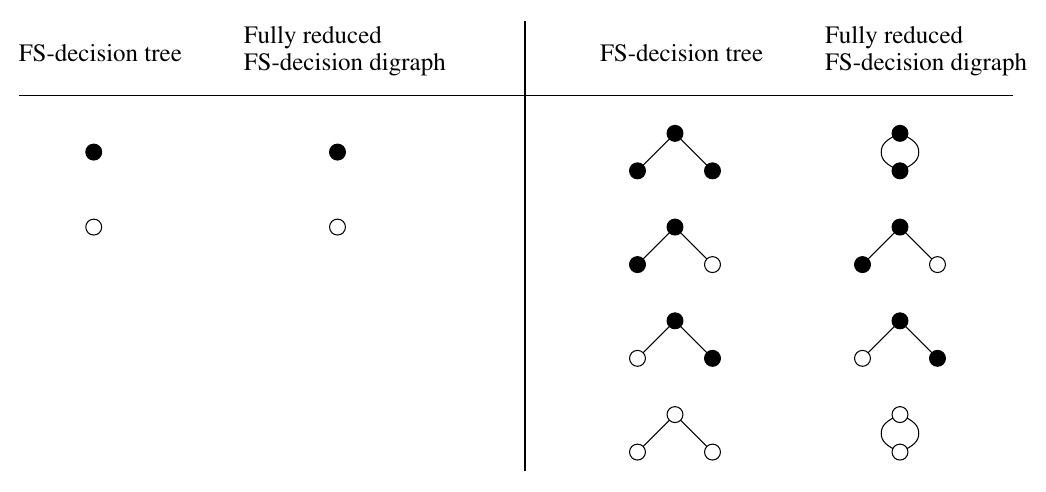}}
\caption{FS-decision trees for binary questionnaires with at most two questions,  and the corresponding fully reduced FS-decision digraphs. All edges are directed downwards, out-neighbours are ordered from left to right, and flagged and unflagged vertices are coloured white and black, respectively. }\label{fig:level-1-2}
\end{figure}

\begin{figure}[t!]
\centerline{\includegraphics[scale=0.7]{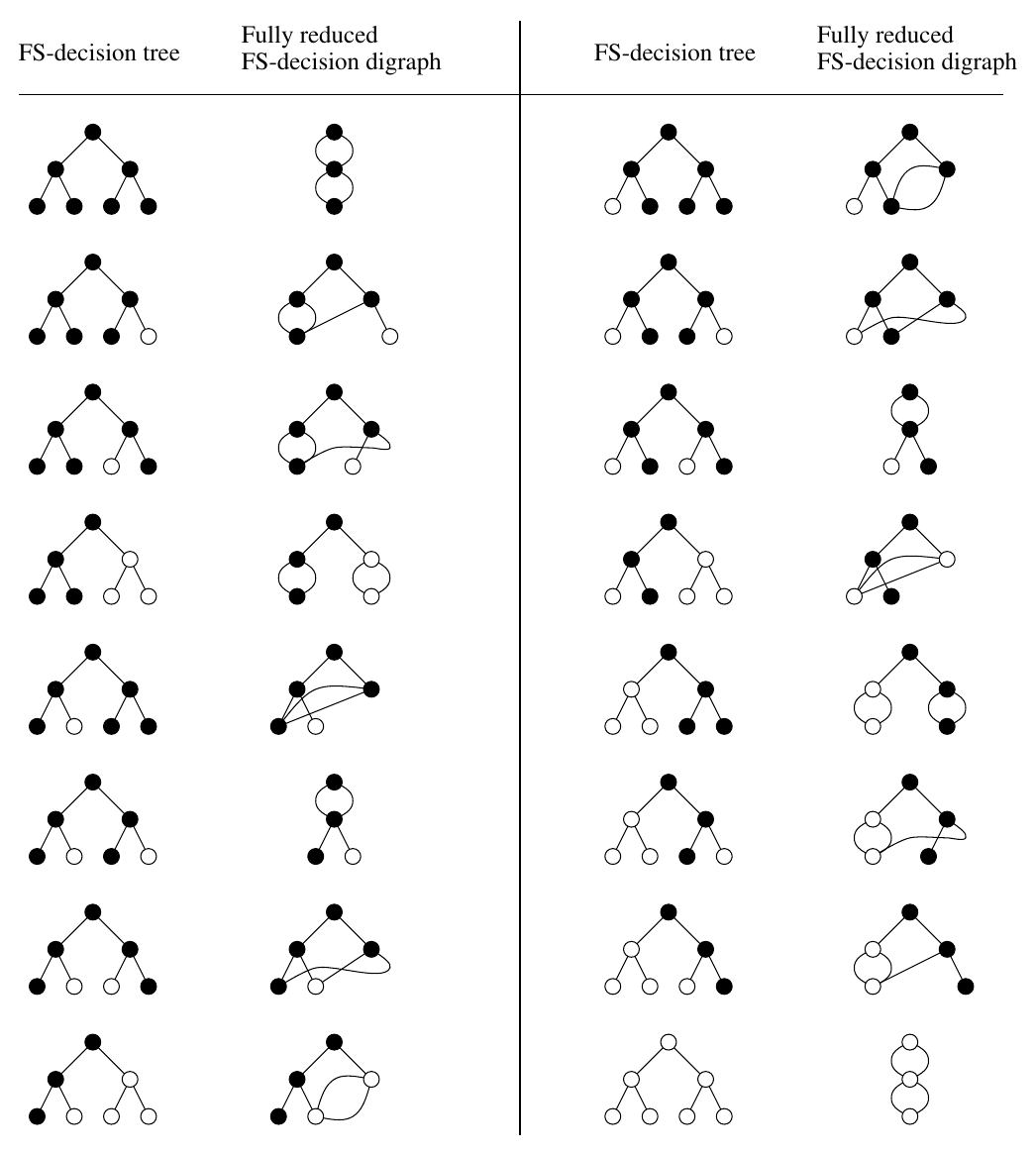}}
\caption{FS-decision trees for binary questionnaires with three questions,  and the corresponding fully reduced FS-decision digraphs. All edges are directed downwards, out-neighbours are ordered from left to right, and flagged and unflagged vertices are coloured white and black, respectively. Continued in Figure~\ref{fig:level-3b}.}\label{fig:level-3}
\end{figure}

\begin{figure}[t!]
\centerline{\includegraphics[scale=0.7]{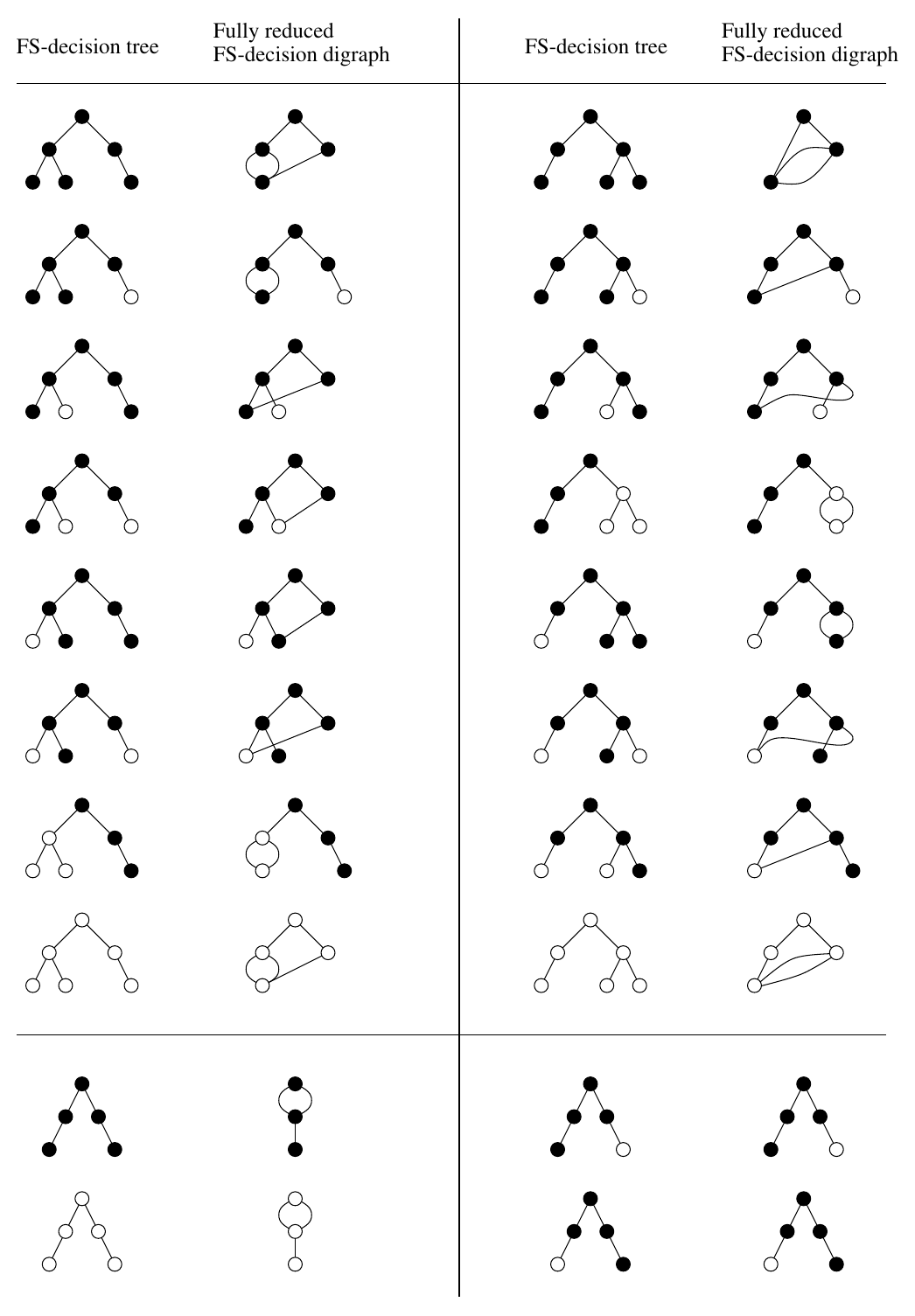}}
\caption{FS-decision trees for binary questionnaires with three questions,  and the corresponding fully reduced FS-decision digraphs. All edges are directed downwards, out-neighbours are ordered from left to right, and flagged and unflagged vertices are coloured white and black, respectively. Continued from Figure~\ref{fig:level-3}.}\label{fig:level-3b}
\end{figure}

Furthermore, Figures~\ref{fig:level-1-2}, \ref{fig:level-3}, and \ref{fig:level-3b} give a list of FS-decision trees for all binary questionnaires with at most 3 questions, and their fully reduced FS-decision digraphs.

The obvious way to construct an FS-decision digraph would be as follows. First, construct the FS-decision tree, then establish equivalence between all (relevant) pairs of vertices. Finally, perform a sequence of merges until no more pairs of equivalent vertices remain. In order to save time and space, however, instead of constructing the FS-decision tree first, we will directly construct the FS-decision digraph that would result from merging all pairs of equivalent vertices; that is, we will be able to recognize such pairs of vertices without first constructing the subtrees rooted at them. We will then prove (Theorem~\ref{the:equivalence}) that this procedure results in the fully reduced FS-decision digraph.

\subsubsection{Recognizing vertex equivalence}

We shall now discuss how to determine whether two vertices in an FS-decision tree are equivalent. We would like to do this as efficiently as possible; that is, without prior construction of the subtrees rooted at the two vertices. For that purpose, we introduce the concepts of a local flag-set and local skip-list.

\begin{defn}\label{def:local-eq}{\rm
Let $k \in \NN$, and let $a=a_0 \ldots a_{k}$ be a $(0,k)$-answer string of an abstract questionnaire $(N,\M)$ with flag-set $F$ and skip-list $\S=(S_0,\ldots,S_{N})$.
\begin{enumerate}[(i)]
\item The {\em local flag-set at $a$} is a set $F^a$ of $(k+1,N-1)$-answer strings defined as follows: if for some $i \le k$, the $(0,N-1)$-answer string $a_0 \ldots a_i* \ldots * \in F$, then $F^a =\{* \ldots * \}$; otherwise,
    $$F^a=\{ b \in A_{k+1}^* \times \ldots \times A_{N-1}^*: ab \in F \}.$$
\item The {\em local skip-list at $a$} is a list $\S^a=(S_{k+1}^a,\ldots,S_{N}^a)$ such that for each $i=k+1,\ldots,N$,
    $$S_{i}^a=\{ b \in A_{k+1}^* \times \ldots \times A_{i-1}^*: ab \in S_{i}  \}.$$
    Note that if $a \in S_{k+1}$, then $S_{k+1}^a=\{ \epsilon \}$; otherwise, $S_{k+1}^a=\emptyset$.
\end{enumerate}
}
\end{defn}

In Theorem~\ref{the:equivalence} below, which will be crucial for constructing FS-decision digraphs, we show that two vertices in an FS-decision tree are equivalent if and only if their local flag-sets, as well as their local skip-lists, are identical. However, for this to hold, the flag-set of the abstract questionnaire must be compatible with its skip-list, in the sense of Definition~\ref{def:compatibleF}.

\begin{thm}\label{the:equivalence}
Let $\Q=(N,\M)$ be an abstract questionnaire  with a flag-set $F$ and skip-list $\S=(S_0,\ldots,S_{N})$. Assume that $m_q \ge 2$ for all $q=0,\ldots,N-1$, and that $F$ is compatible with $\S$. Let $T$ be the FS-decision tree for $\Q$ with question assignment $\kappa$ and answer string assignment $\alpha$, and let $v,w \in V(T)$ be such that $\kappa(v)=\kappa(w)=k+1$. Furthermore, let $\alpha(v)=a=a_0 \ldots a_k$ and $\alpha(w)=b=b_0 \ldots b_k$.

Then $v$ and $w$ are equivalent in $T$ if and only if $F^a=F^b$ and $\S^a=\S^b$.
\end{thm}

\begin{proof}
($\Rightarrow$)  Assume $v$ and $w$ are equivalent in $T$. Thus, by Lemma~\ref{lem:phi}, there exists an isomorphism $\Phi: T_v \to T_w$ satisfying Properties (a) and (b) in Definition~\ref{def:eq-vx}.

First, we show that $F^a=F^b$. Suppose first that $F^a=\{ * \ldots * \}$. Then, by Definition~\ref{def:local-eq}, for some $i \le k$, the $(0,N-1)$-answer string $a_0 \ldots a_i * \ldots * \in F$. Hence $v$ is flagged, and since $w=\Phi(v)$, so is $w$. That means that for some $j \le k$, the $(0,N-1)$-answer string $b_0 \ldots b_j * \ldots * \in F$, and hence $F^b=\{ * \ldots *\}=F^a$.

Hence assume that $F^a \ne \{ * \ldots *\}$. This means that
there is no $i \le k$ such that the $(0,N-1)$-answer string $a_0 \ldots a_i * \ldots * \in F$, and consequently,
the answer string $a_0 \ldots a_k$ is not flagged. As $a_0 \ldots a_k = \alpha(v)$ and $w=\Phi(v)$, we know that $w$ and $\alpha(w)$ are  not flagged. Hence there is no $i \le k$ such that the $(0,N-1)$-answer string $b_0 \ldots b_i * \ldots * \in F$.

Take any $a_{k+1} \ldots a_{N-1} \in F^a$. Then $a_0 \ldots a_k a_{k+1} \ldots a_{N-1} \in F$. Let $\ell \in \{ 0, \ldots N-1 \}$ be such that
$a_0 \ldots a_{N-1}=a_0 \ldots a_{\ell} * \ldots *$ and $\ell$ satisfies the conditions in Definition~\ref{def:compatibleF}. Note that, by the above observation, $\ell >k$. By Lemma~\ref{lem:flagged}, there exists a vertex $u \in V(T)$ such that $\alpha(u)=a_0 \ldots a_{\ell}$. Clearly $u \in V(T_v)$ and $u$ is flagged. Let $z=\Phi(u)$. Then $\alpha(z)=b_0 \ldots b_k a_{k+1} \ldots a_{\ell}$ and $z$ is flagged. Hence for some $j \le \ell$,  the $(0,N-1)$-answer string $b_0 \ldots b_k a_{k+1} \ldots a_j * \ldots * \in F$. Note that by the conclusion of the previous paragraph, we indeed have $j \ge k+1$. By Lemma~\ref{lem:flagged}, there exists a vertex $z' \in V(T)$ such that $\alpha(z')=b_0 \ldots b_k a_{k+1} \ldots a_j$, and necessarily $z' \in V(T_w)$. Let $u'=\Phi^{-1}(z')$. Then $\alpha(u')=a_0 \ldots a_k a_{k+1} \ldots a_j$, and since $z'$ is flagged, so is $u'$. But then the $(0,N-1)$-answer string $a_0 \ldots a_k a_{k+1} \ldots a_i * \ldots * \in F$ for some $k+1 \le i \le j$, and by Definition~\ref{def:compatibleF}(iv), we must have $i \ge \ell$. Since $j \le \ell$, we can see that $i=j=\ell$.

It follows that the $(0,N-1)$-answer string $b_0 \ldots b_k a_{k+1} \ldots a_{\ell} * \ldots * \in F$, and hence the $(k+1,N-1)$-answer string $a_{k+1} \ldots a_{\ell} * \ldots * \in F^b$. Hence $a_{k+1} \ldots a_{N-1} \in F^b$. By symmetry, we conclude that $F^a=F^b$.

Next, we show that $\S^a=\S^b$. If $S_{k+1}^a=\emptyset$, then $a \not\in S_{k+1}$ and hence $v \not\in U$. By Lemma~\ref{lem:phi}(v) we have that $w \not\in U$, and hence $b \not\in S_{k+1}$ and $S_{k+1}^b=\emptyset$. By symmetry, we conclude that $S_{k+1}^a=S_{k+1}^b$. Now let $i \in \{ k+1, \ldots, N-1 \}$, and take any $a_{k+1} \ldots a_{i-1} \in S^a_{i}$. Then $a_0 \ldots a_k a_{k+1} \ldots a_{i-1} \in S_{i}$, and by Lemma~\ref{lem:skip-list-and-vx}, there exists a vertex $u \in V(T)$ such that $\alpha(u)=a_0 \ldots a_k a_{k+1} \ldots a_{i-1}$. Note that $u \in V(T_v)$. Let $z=\Phi(u)$. Then $\alpha(z)=b_0 \ldots b_k a_{k+1} \ldots a_{i-1}$. Since $\alpha(u) \in S_{i}$, we know $u \in U$ and $u$ has a single child. Hence $z$ has a single child, and since $m_{i} \ge 2$, we have $z \in U$. It follows that $\alpha(z) \in S_{i}$, and hence $a_{k+1} \ldots a_{i-1} \in S^b_{i}$. Using symmetry, we obtain $S_i^a=S_i^b$, and thus conclude that  $\S^a=\S^b$.

($\Leftarrow$) Assume $F^a=F^b$ and $\S^a=\S^b$. Define a function
$$\Theta: \{ a \} \cup \left(\{ a \} \times \bigcup_{i=k+1}^{N-1} A_{k+1}^* \times \ldots \times A_i^* \right) \quad \to  \quad \{ b \} \cup \left( \{ b \} \times \bigcup_{i=k+1}^{N-1} A_{k+1}^* \times \ldots \times A_i^* \right)$$ as follows: $\Theta(a)=b$ and
$$\Theta(a a_{k+1} \ldots a_{i})=ba_{k+1} \ldots a_{i} \quad \mbox{ for all } a_{k+1} \ldots a_{i} \in A_{k+1}^* \times \ldots \times A_i^*, \quad k+1 \le i \le N-1.$$
We show that $\Theta$ induces a bijection from $V(T_v)$ to $V(T_w)$.

Take any $u \in V(T_v)$. If $u=v$, then $\alpha(u)=a$ and $\alpha(w)=b=\Theta(a)=\Theta(\alpha(u))$ by the definition of $\Theta$. Otherwise, $\alpha(u)=a a_{k+1} \ldots a_{\ell}$ for some $a_{k+1} \ldots a_{\ell} \in A_{k+1}^* \times \ldots \times A_{\ell}^*$, $k+1 \le \ell \le N-1$. Let $i \le \ell$ be the largest index such that there exists a vertex $z' \in V(T)$ with $\alpha(z')=b a_{k+1} \ldots a_i$. Clearly, $i \ge k$, and suppose that $i < \ell$. Let $u', u'' \in V(T_v)$ be such that $\alpha(u')=a a_{k+1} \ldots a_i$ and $\alpha(u'')=a a_{k+1} \ldots a_i a_{i+1}$. Note that $a_{i+1}=*$ if and only if $u' \in U$, that is, if and only if $\alpha(u') \in S_{i+1}$, that is, if and only if $a_{k+1} \ldots a_i \in S^a_{i+1}$. However, since there is no vertex $z''$ with $\alpha(z'')=b a_{k+1} \ldots a_i a_{i+1}$, we have that $a_{i+1}=*$ if and only if $z' \not\in U$, that is, if and only if $\alpha(z') \not\in S_{i+1}$, that is, if and only if $a_{k+1} \ldots a_i \not\in S^b_{i+1}$. Since $S^a_{i+1}=S^b_{i+1}$, we have a contradiction. Hence $i=\ell$, which means that there exists a vertex $z \in V(T)$ such that $\alpha(z)=\Theta(\alpha(u))$. Note that necessarily $z \in V(T_w)$.

By symmetry, we conclude that $\Theta$ induces a bijection $\Phi: V(T_v) \to V(T_w)$ that satisfies Requirement (a) from Definition~\ref{def:eq-vx}. It remains to show that every $u \in V(T_v)$ is flagged if and only if $\Phi(u)$ is flagged.

Take any $u \in V(T_v)$, and assume $u$ is flagged. First suppose $u=v$. Then there is an index $i$, $0 \le i \le k$, such that the $(0,N-1)$-answer string $a_0 \ldots a_i * \ldots * \in F$. It follows that $F^a=\{ * \ldots * \}$, and hence $F^b=\{ * \ldots * \}$. Consequently, for some $j$, $0 \le j \le k$, we have that the $(0,N-1)$-answer string $b_0 \ldots b_j * \ldots * \in F$, which implies that $\Phi(u)$ is flagged.

Hence assume $u \ne v$, and so $\alpha(u)= a a_{k+1} \ldots a_{\ell}$ for some $a_{k+1} \ldots a_{\ell} \in A_{k+1}^* \times \ldots \times A_{\ell}^*$, $k+1 \le \ell \le N-1$. Since $u$ is flagged, for some $i$, $0 \le i \le \ell$, we have that the $(0,N-1)$-answer string $a_0 \ldots a_i * \ldots * \in F$. If $i \le k$, then it follows that $F^a=\{ * \ldots * \}$, and hence $F^b=\{ * \ldots * \}$. Consequently, for some $j \le k$, the $(0,N-1)$-answer string $b_0 \ldots b_j * \ldots * \in F$, and it follows that $\Phi(u)$ is flagged.

Hence assume $i \ge k+1$.  Since the $(0,N-1)$-answer string $a_0 \ldots a_i * \ldots * \in F$, we have that
the $(k+1,N-1)$-answer string $a_{k+1} \ldots a_i * \ldots *$ is in $F^a$, and so also in $F^b$. Therefore $b a_{k+1} \ldots a_i * \ldots * \in F$, which implies that $\Phi(u) $ is flagged.

By symmetry, we obtain that Requirement (b) in Definition~\ref{def:eq-vx} holds for $\Phi$ as well. Therefore, vertices $v$ and $w$ are equivalent in $T$.
\end{proof}


\subsubsection{Constructing a fully reduced FS-decision digraph}

Theorem~\ref{the:equivalence} will now be used to construct a fully reduced FS-decision digraph without prior construction of the FS-decision tree, thereby saving time and space.

This FS-decision digraph will be constructed in a BFS order, analogously to the construction of the FS-decision tree $T$ in Algorithm~\ref{alg:S-dt}. When a vertex $u$ is processed, its  out-neighbours (children in $T$) are examined in order. If a potential child $c$ is equivalent to an earlier vertex, say $w$, then $c$ is absorbed into $w$, the answer string for $c$ is adjoined to the set of answer strings for $w$, and $w$ is appended to the list of out-neighbours of vertex $u$. Otherwise, $c$ is created as a new vertex and is placed in the queue to be processed later.

This algorithm is detailed below as Algorithm~\ref{alg:S-dd}. Note that procedure flagged can be found in Algorithm~\ref{alg:S-dt}.

\begin{alg}\label{alg:S-dd}{\rm Constructing a fully reduced FS-decision digraph of an abstract questionnaire

\smallskip

\noindent {\bf procedure} FS-digraph$(N,\M,F,\S)$ \\
\# {\em Input: abstract questionnaire $(N,\M)$ with  flag-set $F$, compatible with the skip-list $\S$.}\\
\# {\em Output: a fully reduced FS-decision digraph $D$ for $(N,\M)$ with the subset of skipped vertices $U$, question assignment $\kappa$, answer string assignment $\A$, and flag function $\phi$.} \\
\# {\em The vertices of $D$ are labelled $0,1,2,\ldots$.} \\
\# {\em {\rm Out}$(u)$ is the (ordered) list of out-neighbours of vertex $u$.} \\
\r $\kappa(0):=0$ \\
\r $\A(0):=\{ \epsilon \}$ \\
\r $\phi(0):=0$  \\
\r $L:=[0]$ \qquad \# {\em BFS queue of unprocessed vertices}\\
\r $c:=0$ \qquad \# {\em last vertex label used}\\
\r {\bf for all} $k \in \ZZ_{N+1}$ {\bf do} first$(k):=-1$ \qquad \# {\em label of first vertex with question $k$}\\
\r first$(0):=0$ \\
\r {\bf while} $L \ne [\;]$ {\bf do} \\
\r\r $u:=$ first vertex in $L$ \\
\r\r remove $u$ from $L$ \\
\r\r $q:=\kappa(u)$ \\
\r\r Out$(u):=[\;]$ \\
\r\r {\bf if} $q <N$ {\bf then} \\
\r\r\r {\bf if} $u \in U$ {\bf then} $A=\{ * \}$\\
\r\r\r {\bf else} $A=\ZZ_{m_q}$\\
\r\r\r {\bf for all} $a_q \in A$ {\bf do} \\
\r\r\r\r $c:=c+1$  \qquad \# {\em create a potential new vertex}\\
\r\r\r\r new$:=$ True \\
\r\r\r\r $\kappa(c)=q+1$ \\
\r\r\r\r $\A(c):=\{b a_q: b \in \A(u) \}$ \qquad \# {\em set of answer strings for vertex $c$} \\
\r\r\r\r {\bf if } first$(q+1) \ne -1$ \\
\r\r\r\r {\bf then} \\
\r\r\r\r\r {\bf if } $\exists w \in \{$first$(q+1),\ldots,c-1 \}$ such that equiv$(c,w)=1$ \\
\r\r\r\r\r {\bf then} \\
\r\r\r\r\r\r \# {\em merge vertex $c$ with vertex $w$} \\
\r\r\r\r\r\r new$:=$ False \\
\r\r\r\r\r\r $\A(w):=\A(w) \cup \A(c)$ \\
\r\r\r\r\r\r $c:=c-1$ \\
\r\r\r\r\r\r Out$(u):=$Out$(u) + [ w ]$ \\
\r\r\r\r {\bf else } first$(q+1):=c$ \\
\r\r\r\r {\bf if } new {\bf then} \\
\r\r\r\r\r $a:=$ the first element of $\A(c)$ \\
\r\r\r\r\r {\bf if} $a\in S_{q+1}$  {\bf then} $U := U \cup \{ c \}$ \qquad \# {\em $c$ is a skipped vertex} \\
\r\r\r\r\r {\bf if } flagged$(a,F)=1$ {\bf then} $\phi(c):=1$ \\
\r\r\r\r\r {\bf else } $\phi(c):=0$  \\
\r\r\r\r\r Out$(u):=$Out$(u) + [ c ]$ \\
\r\r\r\r\r $L:= L +  [c]$ \\
\r {\bf return} $c$, Out, $U$, $\kappa$, $\A$, $\phi$

\bigskip

\noindent {\bf procedure} equiv$(v,w)$ \\
\# {\em Input: vertices $v$ and $w$ of an FS-decision tree $T$.}\\
\# {\em The abstract questionnaire $(N,\M)$, its skip-list $\S$, and a compatible flag-set $F$ are considered global variables.}\\
\# {\em Output: 1 if $v$ and $w$ are equivalent in $T$, and 0 otherwise.} \\
\r Ans$:=0$ \\
\r {\bf if } $\kappa(v)=\kappa(w)$ \\
\r {\bf then} \\
\r\r $k=\kappa(v)-1$ \\
\r\r $a:=\alpha(v)$ \\
\r\r $b:=\alpha(w)$ \\
\r\r $t:=$ the $(k+1,N-1)$-answer string $* \ldots *$ \\
\r\r {\bf if } flagged$(a,F)=1$ {\bf then } $F^a:=\{ t \}$ \\
\r\r {\bf else} $F^a:= \{ d \in A_{k+1}^* \times \ldots \times A_{N-1}^*: a d \in F \}$ \\
\r\r {\bf if } flagged$(b,F)=1$ {\bf then } $F^b:=\{ t \}$ \\
\r\r {\bf else} $F^b:= \{ d \in A_{k+1}^* \times \ldots \times A_{N-1}^*: b d \in F \}$ \\
\r\r {\bf if } $F^a=F^b$ \\
\r\r {\bf then} \\
\r\r\r {\bf if } $a \in S_{k+1}$ {\bf then } $S^a_{k+1}:= \{ \epsilon \}$ \\
\r\r\r {\bf else } $S^a_{k+1}:= \{ \; \}$ \\
\r\r\r {\bf if } $b \in S_{k+1}$ {\bf then } $S^b_{k+1}:= \{ \epsilon \}$ \\
\r\r\r {\bf else } $S^b_{k+1}:= \{ \; \}$ \\
\r\r\r $i:=k+1$ \\
\r\r\r {\bf while } $S^a_{i}=S^b_{i}$ {\bf and } $i<N-1$ {\bf do }\\
\r\r\r\r $i:=i+1$ \\
\r\r\r\r $S^a_i:= \{ d \in A_{k+1}^* \times \ldots \times A_{i-1}^*: a d \in S_i \}$ \\
\r\r\r\r $S^b_i:= \{ d \in A_{k+1}^* \times \ldots \times A_{i-1}^*: b d \in S_i \}$ \\
\r\r\r {\bf if } $S^a_{i}= S^b_{i}$ {\bf then } Ans$:=1$  \\
\r {\bf return} Ans
}
\end{alg}

\subsection{Generating a skip-list and a compatible flag-set} \label{sec:generate}

In this section, we explain how a questionnaire designer can generate a skip-list and a compatible flag-set for a given questionnaire, starting from an intuitively constructed ``pre-skip-list'' and ``pre-flag-set''. First, we need a concept of a generalized answer string (Definition~\ref{def:pre-skip-list} below).

As before, we assume that we have an abstract questionnaire $\Q=(N,\M)$, and that $A_i=\ZZ_{m_i}$ is the set of possible answers to question $i$. In addition, for all $i \in \ZZ_N$, we define $A_i^\bl=A_i \cup \{ \bl \}$. We think of the symbol $\bl$ as representing ``anything''. (Note that $\bl$ has a similar meaning as $*$, but without the restrictions imposed by Definitions~\ref{def:skip-list} and \ref{def:compatibleF}.)

\begin{defn}{\rm
Let $k$ and $\ell$ be integers, $0 \le k \le \ell \le N-1$.
A {\em generalized $(k,\ell)$-answer string} for an abstract questionnaire $\Q=(N,\M)$ is an element of $A_k^\bl \times A_{k+1}^\bl \times \ldots \times A_{\ell}^\bl$; that is, a string of the form $a_k a_{k+1} \ldots a_{\ell}$, where $a_i \in A_i^\bl$ for all $i=k,k+1,\ldots,\ell$.
}
\end{defn}

A generalized $(k,\ell)$-answer string  $a_k a_{k+1} \ldots a_{i-1} \bl \, a_{i+1} \ldots a_{\ell}$ can be thought of as representing all $(k,\ell)$-answer strings of the form $a_k a_{k+1} \ldots a_{i-1} x a_{i+1} \ldots a_{\ell}$, for $x \in A_i$.

We next explain what we mean by a pre-skip-list.

\begin{defn}\label{def:pre-skip-list}{\rm
A {\em pre-skip-list} for the abstract questionnaire $\Q=(N,\M)$ is an $(N+1)$-tuple $\P=(P_0,P_1,\ldots,P_{N})$ such that the following hold.
\begin{enumerate}[(i)]
\item $P_0 =\emptyset=P_N$.
\item For each $q$, $0 \le q \le N$, we have that $P_q$ is a set of generalized $(0,q-1)$-answer strings.
\end{enumerate}
}
\end{defn}

For each $q$, we think of the set $P_q$ as containing all $(0,q-1)$-answer strings that necessitate the skipping of question $q$. However, to simplify the work of the questionnaire designer, these answer strings can be provided in the form of generalized answer strings, whereby only the relevant answers are specified, and all irrelevant answers are replaced by the symbol $\bl$. For example, if we have an abstract questionnaire $\Q=(N,\M)$ with $N=5$ and $\M=(2,2,3,2,3)$, and we wish question 4 to be skipped when the answers to questions 0 and 2 are 1 and 0, respectively, then we set $P_4=\{ 1 \bl 0 \bl \}$, which represents the set of answer strings $\{ 1 a_1 0 a_3 : a_i \in A_i, i=1,3 \}$.

We next describe a procedure that converts a pre-skip-list to a skip-list that provides equivalent information but satisfies the requirements of Definition~\ref{def:skip-list}, so it can be used in our main algorithms. Intuitively speaking, for each $q=0, \ldots, N$, each generalized answer string in $S_q$ must be expanded into a corresponding set of answer strings that also satisfies Condition (iii) in Definition~\ref{def:skip-list}; that is, this set must contain only those corresponding answer strings that represent actual vertices of the FS-decision tree.

\begin{alg}\label{alg:skip-list}{\rm Constructing a skip-list of an abstract questionnaire from a pre-skip-list

\smallskip

\noindent {\bf procedure} skip-list$(N,\M,\P)$ \\
\# {\em Input: abstract questionnaire $(N,\M)$ with a pre-skip-list $\P$.}\\
\# {\em Output: skip-list $\S$.} \\
\r $S_0:=\emptyset$ \\
\r $S_N:=\emptyset$ \\
\r $S_1:=P_1$ \qquad \# {\em we may assume $\bl \not\in P_1$} \\ 
\r {\bf for} $q:=2$ {\bf to} $N-1$ \\
\r\r $S_q:=\emptyset$ \\
\r\r {\bf for all} $p_0 \ldots p_{q-1} \in P_q$ {\bf do} \\
\r\r\r {\bf if} $p_0=\bl$ {\bf then} $S:=A_0$ \\
\r\r\r {\bf else} $S:=\{ p_0 \}$ \qquad \# {\em $S$ contains answer strings arising from the generalized answer string $p_0 \ldots p_{q-1}$; at step $i$, it contains $(0,i-1)$-answer strings} \\
\r\r\r {\bf for} $i:=1$ {\bf to} $q-1$ {\bf do} \\
\r\r\r\r {\bf if} $p_i=\bl$ {\bf then} \\
\r\r\r\r\r $S':= \{ ts_i: t \in S-S_i, s_i \in A_i \} \cup \{ t*: t \in S \cap S_i \}$ \\
\r\r\r\r {\bf else} $S':= \{ tp_i: t \in S-S_i \}$ \\
\r\r\r\r $S:=S'$ \\
\r\r\r $S_q:=S_q \cup S$ \\
\r {\bf return} $\S=(S_0, \ldots,S_N)$
}
\end{alg}

Next, we turn our attention to flag-sets. First we need the concept of a pre-flag-set, which we define as follows.

\begin{defn}\label{def:pre-flag-set}{\rm
A {\em pre-flag-set} for the abstract questionnaire $\Q=(N,\M)$ is a set of generalized $(0,N-1)$-answer strings, usually denoted by $P$.}
\end{defn}

We think of a pre-flag-set $P$ as containing every $(0,N-1)$-answer string that is contradictory or of special interest to the questionnaire designer. Again, for simplicity, the answer strings in $P$ are given in the form of generalized answer strings. To illustrate, consider again the example where $\Q=(N,\M)$ with $N=5$ and $\M=(2,2,3,2,3)$. Say that, if chosen together in a response, the answers 1 and 2 to questions 0 and 2, respectively, and the answers 0 and 1 to questions 3 and 4, respectively, are of special interest. Then our pre-flag-set is $P=\{1\bl 2\bl\bl, \bl\bl\bl 01\}$, and it represents the set of answer strings $\{1a_12a_3a_4:a_i\in A_i, i=1,3,4\}   \cup \{a_0a_1a_201:a_i\in A_i, i=0,1,2 \}$.

Given a skip-list $\S$ and pre-flag-set $P$, Algorithm~\ref{alg:flag-set} below creates a flag-set $F$ that is compatible with $\S$, and therefore can be used in our main algorithms. The idea is to expand each generalized answer string in $P$ into a set of answer strings containing equivalent information while also satisfying the conditions of compatibility in Definition~\ref{def:compatibleF}.

\begin{alg}\label{alg:flag-set}{\rm Constructing a compatible flag-set of an abstract questionnaire from a pre-flag-set

\smallskip

\noindent {\bf procedure} flag-set$(N,\M,\S,P)$ \\
\# {\em Input: abstract questionnaire $(N,\M)$ with skip-list $\S$ and pre-flag-set $P$.}\\
\# {\em Output: flag-set $F$.} \\
\r $F:=\emptyset$\\
\r {\bf for all} $p_0 \ldots p_{N-1} \in P$ {\bf do}\\
\r\r {\bf if} $p_0=\bl$ {\bf then} $G:=A_0$\\
\r\r {\bf else} $G:=\{p_0\}$ \qquad \# {\em $G$ contains answer strings arising from the generalized answer string $p_0 \ldots p_{N-1}$; at step $i$, it contains $(0,i-1)$-answer strings} \\
\r\r {\bf for all} $g \in G$ {\bf do} \\
\r\r\r {\bf if} is-contained$(0,g,F)$ {\bf then} $G:=G-\{ g \}$ \\
\r\r {\bf for} $i:=1$ {\bf to} $N-1$ {\bf do}\\
\r\r\r {\bf if} $p_i \ldots p_{N-1}=\bl \ldots \bl$ \\
\r\r\r {\bf then} \\
\r\r\r\r $G':= \{ g * \ldots *: g \in G \}$ \qquad \# {\em $G'$ is a set of $(0,N-1)$-answer strings} \\
\r\r\r\r {\bf for all} $g \in G'$ {\bf do} \\
\r\r\r\r\r {\bf if} is-contained$(N-1,g,F)$ {\bf then} $G':=G'-\{ g \}$ \\
\r\r\r\r $G:=G'$ \\
\r\r\r\r {\bf break} \qquad \# {\em exit the for loop} \\
\r\r\r {\bf else} \\
\r\r\r\r {\bf if} $p_i=\bl$ {\bf then} \\
\r\r\r\r\r $G':=\{ga_i : g\in G-S_i, a_i\in A_i\}\cup\{g*:g\in G\cap S_i\}$ \\
\r\r\r\r {\bf else} $G':=\{gp_i:g \in G-S_i\}$\\
\r\r\r\r {\bf for all} $g \in G'$ {\bf do} \\
\r\r\r\r\r {\bf if} is-contained$(i,g,F)$ {\bf then} $G':=G'-\{ g \}$ \\
\r\r\r\r $G:= G'$\\
\r\r {\bf for all} $g \in G$ {\bf do} \\
\r\r\r {\bf for all} $f \in F$ {\bf do} \\
\r\r\r\r {\bf if} replaces$(g,f)$ {\bf then} $F:=F-\{ f \}$ \\
\r\r $F:=F\cup G$\\
\r {\bf return} $F$

\medskip

\noindent {\bf procedure} is-contained$(k,g,F)$ \\
\# {\em Input: $(0,k)$-answer string $g$ and set $F$ of $(0,N-1)$-answer strings.}\\
\# {\em Output: {\rm True} iff the $(0,N-1)$-answer string $g* \ldots *$ is in $F$.} \\
\r {\bf for} $i:=k+1$ {\bf to} $N-1$ {\bf do} $g:=g*$ \\
\r {\bf if} $g \in F$ {\bf then return} True \\
\r {\bf else} {\bf  return} False

\medskip

\noindent {\bf procedure} replaces$(g,f)$ \\
\# {\em Input: $(0,N-1)$-answer strings $g=g_0 \ldots g_{N-1}$ and $f=f_0 \ldots f_{N-1}$.}\\
\# {\em Output: {\rm True} iff $g=g_0 \ldots g_k* \ldots *$ and $f=g_0 \ldots g_{\ell}* \ldots *$ with $k<\ell$.} \\
\r $i:=0$ \\
\r {\bf while} $g_i=f_i$ and $i<N$ {\bf do} $i:=i+1$ \\
\r {\bf if} $i<N$ and $g_i \ldots g_{N-1}=* \ldots *$ {\bf then return} True \\
\r {\bf else} {\bf  return} False
}
\end{alg}

\section{Example}\label{sec:example}

In this section, we give an example of a simple concrete questionnaire, and use it to illustrate the concepts and algorithms presented in this paper.

Consider the following short questionnaire directed at university students taking a certain course.

\begin{framed}

{\footnotesize

\begin{enumerate}
\item[0.] For your program, this course is
\begin{enumerate}
    \item[0)] compulsory
    \item[1)] not required
\end{enumerate}
\item[1.] Which of the following best describes your reason for taking this course?
\begin{enumerate}
    \item [0)] It seemed interesting
    \item[1)] I wanted to challenge myself
    \item[2)] I wanted to boost my GPA
\end{enumerate}
\item[2.] What percentage of classes did you attend?
\begin{enumerate}
    \item[0)] 0-25\%
    \item[1)] 26-50\%
    \item[2)] 51-75\%
    \item[3)] 76-100\%
\end{enumerate}
\item[3.] I think the professor conveys the subject matter effectively.
\begin{enumerate}
    \item [0)] agree
    \item[1)] disagree
\end{enumerate}
\item[4.] This course is challenging.
\begin{enumerate}
    \item [0)] agree
    \item[1)] disagree
\end{enumerate}
\end{enumerate}
}

\end{framed}

Note that the corresponding abstract questionnaire has $N=5$ and $\M=(2,3,4,2,2)$.

\begin{figure}[t!]
\centerline{\includegraphics[scale=0.7]{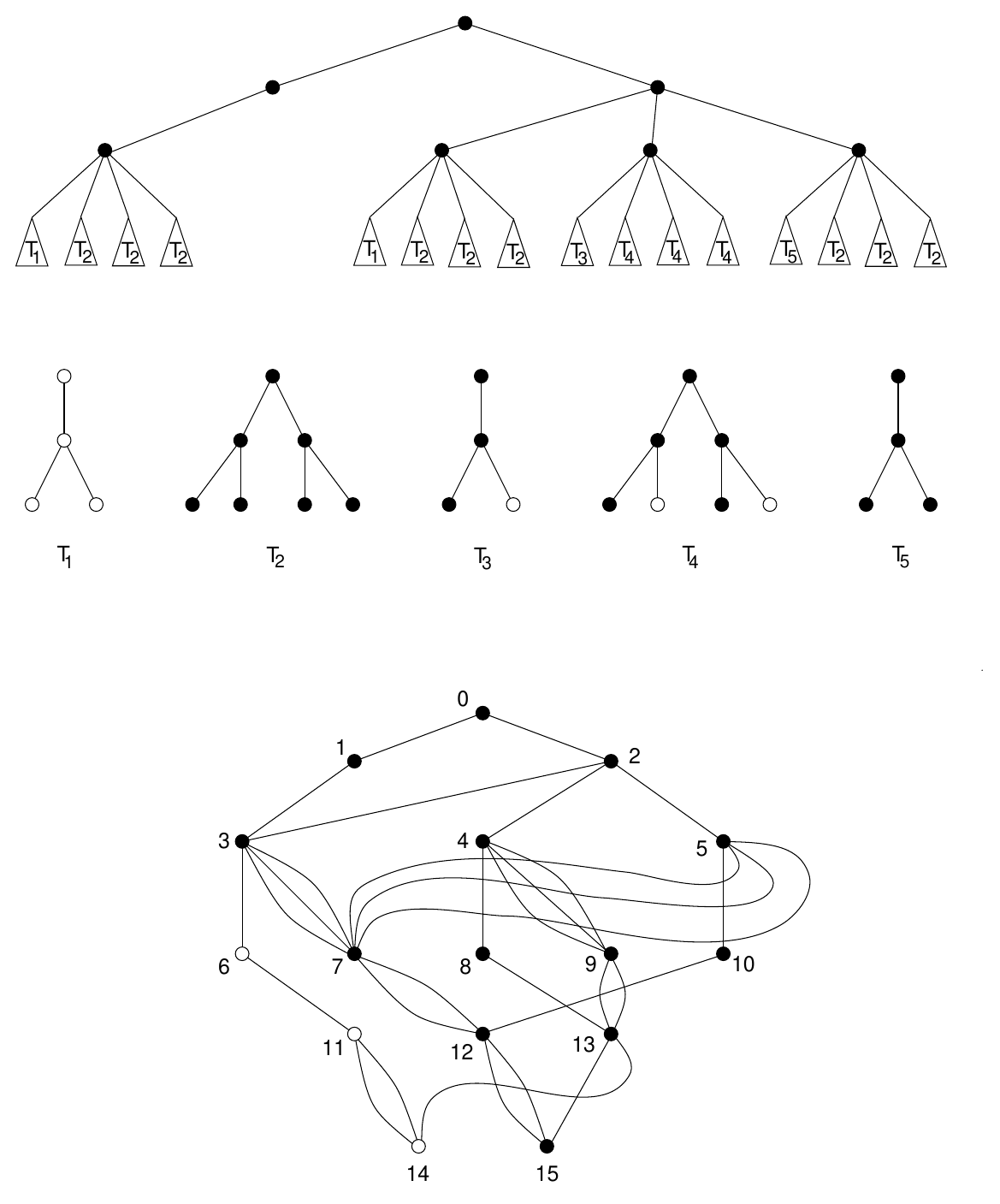}}
\caption{The FS-decision tree (top) and the fully reduced FS-decision digraph (bottom) for the example in Section~\ref{sec:example}. All edges are directed downwards, out-neighbours are ordered from left to right, and flagged and unflagged vertices are coloured white and black, respectively. Vertices in the FS-decision digraph are labelled in the order created by Algorithm~\ref{alg:S-dd}.}\label{fig:grand-eg}
\end{figure}

First, we shall illustrate the results of Section~\ref{sec:ordering}. Suppose that we insist that questions 1 and 2 be asked after question 0, and that questions 4 and 3 be asked after questions 1 and 2, respectively. This gives us the precedence relation $R=\{(0,1),(0,2),(1,4),(2,3)\}$. Note that $R$ is an irreflexive binary relation with no directed cycles. By Lemma~\ref{lem:total}, there exists a total order extending $R$, and we may apply Algorithm~\ref{alg:Q}. The resulting output is
$$\T=[\;[0,1,2,3,4], [0,1,2,4,3],[0,1,4,2,3],[0,2,1,3,4],
[0,2,1,4,3], [0,2,3,1,4]\;];$$
that is, $\T$ is a list of all total orders on $\ZZ_5$ that extend the relation $R$.
Note that for the remainder of this section, we will be using the default question order $[0,1,2,3,4]$.

Next, we would like to represent the flow of our questionnaire with a graph. To illustrate the work of Section~\ref{sec:main}, suppose that we would like to impose the following restrictions on our questionnaire: if a respondent answers with 0 to question 0, they should skip question 1 (no need to ask why they took the course since it was compulsory for them), and if they answer with 0 to question 2, they should skip question 3 (since they attended very few classes, their opinion on the professor is not very relevant). Thus our pre-skip-list is $\P=(P_0,P_1, \ldots,P_5)$, where $P_1=\{0\}$, $P_3=\{\bl\bl0\}$, and $P_0=P_2=P_4=P_5=\emptyset$. Applying Algorithm~\ref{alg:skip-list} gives us the skip list $\S=(S_0,S_1,\ldots,S_5)$, where $S_1=\{0\}$, $S_3=\{0*0,100,110,120\}$, and $S_0=S_2=S_4=S_5=\emptyset$.

We would also like to flag certain answer strings. It is questionable for someone to take the course out of interest, then not attend any classes, so we add $\bl 0 0 \bl \bl$ to the pre-flag-set $P$. We would also like to know more about why someone who wanted to challenge themselves with this course did not find the course challenging; hence, we let $\bl 1 \bl \bl 1 \in P$. Lastly, we would like to keep track of those who are required to take the course yet who have attended very few classes, so we add $0 \bl 0 \bl \bl$ to the pre-flag-set as well. Thus, we have a pre-flag-set $P=\{\bl00\bl\bl,\bl1\bl\,\bl\,1,0\bl0\bl\bl\}$. To create a flag-set $F$ that is compatible with the skip-list $\S$, we expand each of the answer strings as outlined in Algorithm~\ref{alg:flag-set}. We thus obtain $$F=\{ 100**,110*1,11101,11111,11201,11211,11301,11311,0*0**\}.$$

Finally, we apply Algorithm~\ref{alg:S-dt} to obtain the FS-decision tree for $\Q$, and Algorithm~\ref{alg:S-dd} to obtain the fully reduced FS-decision digraph for $\Q$ (see Figure~\ref{fig:grand-eg}).

\section{Conclusion}

In this paper, we introduced FS-decision trees and, more importantly, FS-decision digraphs as new models for abstract questionnaires endowed with the additional information in the form of a skip-list and a flag-set. We presented algorithms for constructing both models, as well as for generating suitable input data (a skip-list and a compatible flag-set) from the more intuitive pre-skip-list and pre-flag-set. We also described how to construct all possible orderings of the questions based solely on an abstract precedence relation. Our hope is that our models will help questionnaire designers visualize and automatize their work.

\bigskip\bigskip

\centerline{\bf Acknowledgement}

\medskip

The third author gratefully acknowledges support by the Natural Sciences and Engineering Research Council of Canada (NSERC), Discovery Grant RGPIN-2022-02994.

\end{document}